\documentclass[12pt, a4paper, twoside]{article}
\usepackage{amssymb}
\usepackage[leqno]{amsmath}
\usepackage{latexsym}
\usepackage{enumerate}
\usepackage{amsmath, amssymb, amscd}
\Roman{equation}

\begin{document}

\title{\textbf{\Large{The classification of ACM line bundles on quartic hypersurfaces on $\mathbb{P}^3$}}}

\author{Kenta Watanabe \thanks{Department of Mathematical Sciences, Osaka University, 1-1 Machikaneyama-chou Toyonaka Osaka 560-0043 Japan, {\it E-mail address:goo314kenta@mail.goo.ne.jp}, Telephone numbers: 090-9777-1974} }

\date{}

\maketitle


\begin{abstract}

In this paper, we give a complete classification of initialized and ACM line bundles on a smooth quartic hypersurface on $\mathbb{P}^3$.

\end{abstract}

\section{Introduction}

We work over the complex number field $\mathbb{C}$. Let $X$ be a smooth hypersurface on $\mathbb{P}^n$ and $\mathcal{O}_X(1)$ be the very ample line bundle given by a hyperplane section of $X$. Then, a vector bundle $\mathcal{E}$ on $X$ is called an {\it{Arithmetically Cohen-Macaulay}} ({\it{ACM}} for short) if the cohomology groups $H^i(X,\mathcal{E}(l))$ for $1\leq i\leq\dim(X)-1$ and $l\in\mathbb{Z}$ vanish, where $\mathcal{E}(l)=\mathcal{E}\otimes\mathcal{O}_X(l)$. It is well known that if a vector bundle $\mathcal{E}$ on $X$ splits, that is, $\mathcal{E}$ is a direct sum of line bundles, then $\mathcal{E}$ is an ACM bundle. If $X=\mathbb{P}^n\;(n\geq1)$, the converse assertion of it is also correct. However, in the case where $X$ is a hypersurface of degree $d\geq2$ on $\mathbb{P}^n$, there exists a counterexample. Therefore, for an ACM bundle on hypersurfaces of higher degree, we can consider the splitting problem.

In the past decades, several authors have studied about indecomposable ACM bundles on a hypersurface on $\mathbb{P}^n$. For example, Kn${\rm{\ddot{o}}}$rrer [Kn] has completely classified ACM bundles on quadrics. In [C-H], Casanellas and Hartshorne have constructed a $n^2+1$- dimensional family of rank $n$ indecomposable ACM vector bundles with Charn classes $c_1=nH$ and $c_2=\frac{1}{2}(3n^2-n)$ on a smooth cubic surface in $\mathbb{P}^3$ for $n\geq2$. Moreover, Faenzi [Fa] gave a precise classification of rank 2 ACM bundles on cubic surfaces. 

In this article, we are interested in ACM bundles on a hypersurface $X$ of degree 4 on $\mathbb{P}^3$. In order to investigate such bundles of higher rank, we have classified ACM and initialized line bundles (i.e., line bundles $\mathcal{L}$ with $H^0(X,\mathcal{L})\neq0$ and $H^0(X,\mathcal{L}(-1))=0$). Our main theorem is as follows.

\newtheorem{thm}{Theorem}[section]

\begin{thm} Let $X$ be a smooth quartic hypersurface on $\mathbb{P}^3$, and let $D$ be a non-zero effective divisor on $X$. Then the following conditions are equivalent.

\smallskip

\smallskip

\noindent {\rm{(i)}} $\mathcal{O}_X(D)$ is an ACM and initialized line bundle.

\noindent {\rm{(ii)}} For any general member $C\in |\mathcal{O}_X(1)|$, one of the following cases occurs.

\smallskip

{\rm{(a)}} $D^2=-2$ and $1\leq C.D\leq3$.

{\rm{(b)}} $D^2=0$ and $3\leq C.D\leq4$.

{\rm{(c)}} $D^2=2$ and $C.D=5$.

{\rm{(d)}} $D^2=4$, $C.D=6$ and $|D-C|=|2C-D|=\phi$.

\end{thm}

Our plan of this paper is as follows. In Section 2, we recall some classical results about linear systems and ACM line bundles on K3 surfaces, and prepare some lemmas to prove the main theorem. In Section 3, we give a proof of Theorem 1.1.

\smallskip

\smallskip

\noindent {\bf{Notation and Conventions.}} A surface is a smooth projective surface. Let $X$ be a surface. For a divisor $D$ on $X$, we denote by $|D|$ the linear system defined by $D$. If two divisors $D_1$ and $D_2$ on $X$ satisfy the condition that $|D_1|=|D_2|$, then we will write $D_1\sim D_2$. We will denote by $\mathcal{O}_X(1)$ a very ample line bundle which provides a closed embedding in some $\mathbb{P}^n$ and denote by $\mathcal{O}_X(l)$ the line bundle $\mathcal{O}_X(1)^{\otimes l}$. For a vector bundle $\mathcal{E}$ on $X$, we will write $\mathcal{E}\otimes\mathcal{O}_X(l)=\mathcal{E}(l)$. We will say that a vector bundle $\mathcal{E}$ on $X$ is {\it{initialized}} with respect to $\mathcal{O}_X(1)$ if it satisfies the following condition.
$$ H^0(X,\mathcal{E}(-1))=0 \text{ but } H^0(X,\mathcal{E})\neq0.$$
\noindent We call a regular surface $X$ a K3 surface if the canonical bundle of $X$ is trivial. Note that a quartic hypersurface on $\mathbb{P}^3$ is a K3 surface.

\section{Linear systems and ACM line bundles on K3 surfaces}

In this section, we recall some classical results about linear systems on K3 surfaces. Moreover, we state an important feature of ACM line bundles on K3 surfaces, and give some lemmas to prove Theorem 1.1.

\newtheorem{prop}{Proposition}[section]

\begin{prop}{\rm{([SD], Lemma 3.7)}} Let $X$ be a K3 surface, and $L$ be a base point free and big line bundle. Then, any member of $|L|$ is 2-connected.\end{prop}

\noindent For a non-zero effective divisor $D$ with $D^2\geq0$ on a K3 surface and the base divisor $\Delta$ of $D$, we have $D-\Delta\neq0$. Therefore, by Proposition 2.1, for such a divisor $D$ and any base point free and big line bundle $L$, we have $D.L\geq2$.

$\;$

\noindent {\bf{Remark 2.1}}. Let $X$ be a K3 surface, and let $\mathcal{L}$ be a line bundle on $X$. Then, $\mathcal{L}$ is ACM if and only if the dual of it is ACM.

\smallskip

\smallskip

\noindent This assertion can be easily proved by using Serre duality.

\newtheorem{lem}{Lemma}[section]

\begin{lem} Let $X$ be a quartic surface on $\mathbb{P}^3$, and let $D$ be a non-zero effective divisor with $D^2\geq0$ on $X$. Then, $D.\mathcal{O}_X(1)\geq3$.\end{lem}

{\it{Proof}}. Let $C\in |\mathcal{O}_X(1)|$ be a smooth curve. Note that $C$ is a plane quartic curve. If $D^2\geq2$, then, by Hodge index theorem, we have
$$(D.C)^2\geq D^2C^2\geq 8.$$
\noindent Hence, the assertion holds. We consider the case where $D^2=0$. Then, by Proposition 2.1, we have $D.C\geq2$. Assume that $D.C=2$. By Proposition 2.1, for the base divisor $\Delta$ of $|D|$, we have $C.(D-\Delta)\geq2$, and hence $C.(D-\Delta)=2$ and $C.\Delta=0$. Since $C$ is ample, $|D|$ is base point free and the general member of $|D|$ is an elliptic curve and hence, we have $h^1(\mathcal{O}_X(D))=0$. 

On the other hand, since $(C-D)^2=0$ and $C.(C-D)=2$, by the same reason, the general member of $|C-D|$ is an elliptic curve. Hence, we have $h^1(\mathcal{O}_X(C-D))=0$. By the exact sequence
$$0\rightarrow\mathcal{O}_X(D-C)\rightarrow\mathcal{O}_X(D)\rightarrow\mathcal{O}_C(D)\rightarrow0,$$
\noindent we have $h^0(\mathcal{O}_C(D))=h^0(\mathcal{O}_X(D))=2$, and hence, $|D|_C|$ is a pencil of degree $C.D=2$. However, this contradicts to the fact that $C$ is a trigonal curve. Hence, the assertion holds.$\hfill\square$

$\;$

\noindent {\bf{Corollary 2.1}} {\it{Let}} $X$ {\it{and}} $D$ {\it{be as in Lemma 2.1. Then, if}} $D.\mathcal{O}_X(1)=3$, $|D|$ {\it{is base point free and the general member of it is an irreducible smooth curve}}.

\smallskip

\smallskip 

{\it{Proof}}. Since, for the base divisor $\Delta$ of $|D|$, $(D-\Delta)^2\geq0$, by Lemma 2.1, we have $(D-\Delta).\mathcal{O}_X(1)\geq3$. By the assumption, we have $(D-\Delta).\mathcal{O}_X(1)=3$ and $\Delta.\mathcal{O}_X(1)=0$. Since $\mathcal{O}_X(1)$ is ample, we have $\Delta=0$. If $D^2>0$, then, by Bertini's theorem, we have the assertion. If $D^2=0$, then there exists an elliptic curve $F$ with $D\sim rF\;(r\geq1)$. By Lemma 2.1, we have $D\sim F$. Therefore, the assertion holds.$\hfill\square$

\begin{lem} Let $X$ be a quartic surface on $\mathbb{P}^3$, and $D$ be an effective divisor with $h^1(\mathcal{O}_X(D))=0$. Then, we have the following results.

\smallskip

\smallskip

\noindent {\rm{(i)}} If $h^1(\mathcal{O}_X(-D)(1))=0$ and $\mathcal{O}_X(1).D\leq3$, then $\mathcal{O}_X(D)$ is an ACM line bundle.

\noindent {\rm{(ii)}} If $h^1(\mathcal{O}_X(-D)(1))=h^1(\mathcal{O}_X(-D)(2))=0$ and $\mathcal{O}_X(1).D\leq7$, then $\mathcal{O}_X(D)$ is an ACM line bundle.\end{lem}

{\it{Proof}}. If $D=0$, then the assertion is clear. Hence, we assume that $D\neq0$. Let $C\in |\mathcal{O}_X(1)|$ be a smooth curve. First of all, for $n\geq0$, we have
$$h^1(\mathcal{O}_C((n+1)C+D))=h^0(\mathcal{O}_C(-nC-D))=0.$$
\noindent By the exact sequence
$$0\rightarrow\mathcal{O}_X(nC+D)\rightarrow\mathcal{O}_X((n+1)C+D)\rightarrow\mathcal{O}_C((n+1)C+D)\rightarrow0,$$
\noindent for any $m\geq0$, we have 
$$h^1(\mathcal{O}_X(-mC-D))=h^1(\mathcal{O}_X(mC+D))=0,$$
\noindent by using mathematical induction. By Remark 2.1, it is sufficient to show that $\mathcal{O}_X(-D)$ is ACM.

\smallskip

\smallskip

(i) For $n\geq1$, by the assumption, we have
$$C(D-nC)\leq3-4n\leq-1.$$
\noindent Hence, we have
$$h^1(\mathcal{O}_C((n+1)C-D))=h^0(\mathcal{O}_C(D-nC))=0.$$
\noindent Since $h^1(\mathcal{O}_X(C-D))=0$, by the exact sequence
$$0\rightarrow\mathcal{O}_X(nC-D)\rightarrow\mathcal{O}_X((n+1)C-D)\rightarrow\mathcal{O}_C((n+1)C-D)\rightarrow0,$$
\noindent for $m\geq1$, we have $h^1(\mathcal{O}_X(mC-D))=0$, by using mathematical induction. Hence, the assertion holds.

\smallskip

\smallskip

(ii) For $n\geq2$, by the assumption, we have
$$C(D-nC)\leq7-4n\leq-1.$$
\noindent Hence, we have
$$h^1(\mathcal{O}_C((n+1)C-D))=h^0(\mathcal{O}_C(D-nC))=0.$$
\noindent Since $h^1(\mathcal{O}_X(C-D))=h^1(\mathcal{O}_X(2C-D))=0$, by the exact sequence
$$0\rightarrow\mathcal{O}_X(nC-D)\rightarrow\mathcal{O}_X((n+1)C-D)\rightarrow\mathcal{O}_C((n+1)C-D)\rightarrow0,$$
\noindent for $m\geq1$, we have $h^1(\mathcal{O}_X(mC-D))=0$, by using mathematical induction. Hence, the assertion holds.$\hfill\square$

\section{Proof of the main result}

In this section, we prove Theorem 1.1. Here, notations are as in Theorem 1.1.

$\;$

{\it{Proof of Theorem 1.1}}. Let $C\in|\mathcal{O}_X(1)|$ be a smooth curve.

\smallskip

\smallskip

\noindent (i)$\Longrightarrow$(ii) Note that since $\mathcal{O}_X(D)$ is initialized, $D\notin|\mathcal{O}_X(1)|$, and that since $\mathcal{O}_X(D)$ is ACM, we have $h^1(\mathcal{O}_X(D))=0$. Since $D$ is effective, $\chi(\mathcal{O}_X(D))\geq1$ and hence, $D^2\geq-2$.

\smallskip

We consider the case where $D^2=-2$. Assume that $|C-D|\neq\phi$. By the ampleness of $C$, we have $C.(C-D)>0$ and hence, we have $1\leq C.D\leq3$. If $C.D=3$, we have $(C-D)^2=-4$. Hence, by Riemann-Roch theorem, we have $h^1(\mathcal{O}_X(C-D))>0$. However, this contradicts to the assumption. Hence, in this case, we have $1\leq C.D\leq 2$. Assume that $|C-D|=\phi$. Since $\mathcal{O}_X(D)$ is initialized and ACM, we have $|D-C|=\phi$ and $h^1(\mathcal{O}_X(C-D))=0$. Hence, by Riemann-Roch theorem, we have $C.D=3$.

\smallskip

We consider the case where $D^2=0$. Assume that $|C-D|\neq\phi$. By the same reason, we have $1\leq C.D\leq3$. Since $D^2=0$, by Lemma 2.1, we have $C.D=3$. Assume that $|C-D|=\phi$. Since $\mathcal{O}_X(D)$ is ACM and initialized, by Riemann-Roch theorem, we have $C.D=4$. 

\smallskip

We consider the case where $D^2=2$. Assume that $|C-D|\neq\phi$. Since we have $1\leq C.D\leq3$, by Lemma 2.1, we have $C.D=3$. Although we have $(C-D)^2=0$ and $C(C-D)=1$, this contradicts to the assertion of Lemma 2.1. Therefore, we have $|C-D|=\phi$. By the assumption, we have $C.D=5$.

\smallskip

We consider the case where $D^2\geq4$. In this case, we have $|C-D|=\phi$. In fact, if we assume that $|C-D|\neq\phi$, then, by the same reason, we have $1\leq C.D\leq3$. However, by Hodge index theorem, we have
$$16\leq C^2.D^2\leq (C.D)^2\leq9.$$
\noindent This is a contradiction. Since $\mathcal{O}_X(D)$ is ACM and initialized, we have $|D-C|=\phi$ and
$$\chi(\mathcal{O}_X(C-D))=0.$$
\noindent Hence, we have 
$$D^2=2C.D-8.$$
\noindent On the other hand, by the same reason, we have
$$\chi(\mathcal{O}_X(2C-D))\geq0.$$
\noindent Therefore, have $C.D\leq6$ and $D^2\leq4$ and hence, we have $C.D=6$ and $D^2=4$. Since $\chi(\mathcal{O}_X(2C-D))=0$, we also have $|2C-D|=\phi$.

$\;$

(ii)$\Longrightarrow$(i) We consider the case where $D^2=-2$. Assume that $C.D=1$. Since $C$ is ample, $D$ is irreducible and hence, we have $h^1(\mathcal{O}_X(D))=0$. Since $(C-D)^2=0$ and $C.(C-D)=3$, by Corollary 2.1, $|C-D|$ is base point free and the general member of it is irreducible. Hence, we have $h^1(\mathcal{O}_X(C-D))=0$. Therefore, by Lemma 2.2, $\mathcal{O}_X(D)$ is ACM. On the other hand, since $h^0(\mathcal{O}_X(D-C))=0$, $\mathcal{O}_X(D)$ is initialized. 

Assume that $C.D=2$. Then the movable part of $|D|$ is empty, otherwise, by the ampleness of $C$ and Lemma 2.1, we have $C.D>3$. Hence, we have $h^0(\mathcal{O}_X(D))=1$. By Riemann-Roch theorem, we have $h^1(\mathcal{O}_X(D))=0$. On the other hand, since we have $(C-D)^2=-2$ and $C.(C-D)=2$, by the same reason, we have $h^1(\mathcal{O}_X(C-D))=0$. Hence, by Lemma 2.2, $\mathcal{O}_X(D)$ is ACM. Since we have $h^0(\mathcal{O}_X(D-C))=0$, $\mathcal{O}_X(D)$ is initialized. 

Assume that $C.D=3$. by the same reason, the movable part of $|D|$ is empty. Hence, we have $h^0(\mathcal{O}_X(D))=1$ and $h^1(\mathcal{O}_X(D))=0$. On the other hand, assume that $|C-D|\neq\phi$. Since $(C-D).(4C-D)=-1$, we have
$$h^0(\mathcal{O}_X(3C))=h^0(\mathcal{O}_X(4C-D)).$$
\noindent Since $h^1(\mathcal{O}_X(3C))=0$, we have
$$20=h^0(\mathcal{O}_X(3C))\geq \chi(\mathcal{O}_X(4C-D))=21.$$
\noindent This is a contradiction. Hence, we have $|C-D|=\phi$. Since $C.(D-C)=-1$, we have $|D-C|=\phi$ and hence, $\mathcal{O}_X(D)$ is initialized. Moreover, we have
$$h^1(\mathcal{O}_X(C-D))=-\chi(\mathcal{O}_X(C-D))=0$$
\noindent and hence, by Lemma 2.2, $\mathcal{O}_X(D)$ is ACM.

\smallskip

We consider the case where $D^2=0$. Assume that $C.D=3$. First of all, Corollary 2.1, the general member of $|D|$ is irreducible. Hence, we have $h^1(\mathcal{O}_X(D))=0$. On the other hand, since $(C-D)^2=-2$ and $C.(C-D)=1$, $|C-D|\neq\phi$ and the member of it is irreducible. Therefore, we have $h^1(\mathcal{O}_X(C-D))=0$. By Lemma 2.2, $\mathcal{O}_X(D)$ is ACM. Since $C.(D-C)=-1$, we have $|D-C|=\phi$ and hence, $\mathcal{O}_X(D)$ is initialized.

\smallskip

Assume that $C.D=4$. First of all, we note that since $(C-D)^2=-4$ and $C(C-D)=C(D-C)=0$, we have $|C-D|=|D-C|=\phi$ and hence, $\mathcal{O}_X(D)$ is initialized and we have
$$h^1(\mathcal{O}_X(C-D))=-\chi(\mathcal{O}_X(C-D))=0.$$

We consider the case where $|D|$ is not base point free. Let $\Delta$ be the base divisor of $|D|$, and let $D^{'}=D-\Delta$. By the ampleness of $C$ and Lemma 2.1, we have
$$3\leq C.D^{'}<C.D=4.$$
\noindent Hence, we have $C.D^{'}=3$ and $C.\Delta=1$ (this implies that $\Delta$ is a $(-2)$-curve). Since $D^2=0$, we have $2D^{'}.\Delta=-{D^{'}}^2+2$. Since $D^{'}.\Delta\geq0$, we have ${D^{'}}^2=0$ or 2. If ${D^{'}}^2=2$, we have $(D^{'}-\Delta)^2=0$. Since $C.(D^{'}-\Delta)=2$, by Lemma 2.1, we have a contradiction. Hence, we have ${D^{'}}^2=0$. Since we have $D^{'}.\Delta=1$, $D$ is 1-connected and hence, we have $h^1(\mathcal{O}_X(D))=0$. Here, we show that $h^1(\mathcal{O}_X(2C-D))=0$. First of all, since $(C-D^{'})^2=-2$ and $C.(C-D^{'})=1$, the member of $|C-D^{'}|$ is irreducible. On the other hand, since $(C-\Delta)^2=0$ and $C.(C-\Delta)=3$, by Corollary 2.1, $|C-\Delta|$ is base point free and the general member of it is irreducible. Since $(C-\Delta).(C-D^{'})=1$, $|2C-D|$ contains a 1-connected divisor and hence, we have $h^1(\mathcal{O}_X(2C-D))=0$. Therefore, by Lemma 2.2, $\mathcal{O}_X(D)$ is ACM. 

We consider the case where $|D|$ is base point free. Then, there exists an elliptic curve $F$ such that $D\sim rF\;(r\geq1)$. Since $C.D=4$, by Lemma 2.1, we have $D\sim F$ and hence, we have $h^1(\mathcal{O}_X(D))=0$. We show that $h^1(\mathcal{O}_X(2C-D))=0$. Assume that $|2C-D|$ is not base point free and let $\Delta$ be the base divisor of it, and let $D^{'}=2C-D-\Delta$. Since $C$ is ample, by Lemma 2.1, we have
$$3\leq D^{'}.C<(2C-D).C=4.$$
\noindent Hence, we have $D^{'}.C=3$ and $\Delta.C=1$ (this implies that $\Delta$ is a $(-2)$-curve). Since $(2C-D)^2=0$, we have $2D^{'}.\Delta=-{D^{'}}^2+2$. Since $D^{'}.\Delta\geq0$, we have ${D^{'}}^2=0$ or 2. If ${D^{'}}^2=2$, then we have $D^{'}.\Delta=0$ and hence, $(C-D^{'}).(C-\Delta)=0$. Since $|C-D^{'}|\neq\phi$ and  $|C-\Delta|\neq\phi$, $D$ is not 1-connected. However, this contradicts to the fact that $h^1(\mathcal{O}_X(D))=0$. If ${D^{'}}^2=0$, we have $D^{'}.\Delta=1$. Then, $|C-D^{'}|\neq\phi$ and $D.(C-D^{'})=-1$. This contradicts to the assumption that $|D|$ is base point free. Therefore, $|2C-D|$ is base point free. Since $C.(2C-D)=4$ and $(2C-D)^2=0$, by the same reason, the general member of $|2C-D|$ is an elliptic curve. Hence, we have $h^1(\mathcal{O}_X(2C-D))=0$. By Lemma 2.2, $\mathcal{O}_X(D)$ is ACM.

\smallskip

We consider the case where $D^2=2$. First of all, we show that $|D|$ is base point free. Assume that $|D|$ is not base point free and let $\Delta$ be the base divisor of $|D|$ and let $D^{'}=D-\Delta$. Since $D^{'}.C<D.C=5$, if ${D^{'}}^2=0$, then by Lemma 2.1, there exists an elliptic curve $F$ with $D^{'}\sim F$. Then, since we have $h^1(\mathcal{O}_X(D^{'}))=0$, we have 
$$\chi(\mathcal{O}_X(D^{'}))=h^0(\mathcal{O}_X(D^{'}))=h^0(\mathcal{O}_X(D))\geq\chi(\mathcal{O}_X(D)).$$
\noindent Therefore, we have ${D^{'}}^2\geq D^2=2$. This is a contradiction. Hence, we have ${D^{'}}^2\geq2$. Since, by Lemma 2.1, we have $C.D^{'}\geq3$, we have $C.D^{'}=3$ or 4. If $C.D^{'}=3$, then we have $(C-D^{'})^2\geq0$ and $C.(C-D^{'})=1$. However, this contradicts to the assertion of Lemma 2.1. Assume that $C.D^{'}=4$. Then, we have
$$C.(C-D^{'})=C.(D^{'}-C)=0.$$
\noindent Since $(C-D^{'})^2\geq-2$, by the ampleness of $C$, we have $C\sim D^{'}$. Since $C.\Delta=1$, $\Delta$ is a $(-2)$-curve. However, since $D^2=2$, we have $C.\Delta=0$. This is a contradiction. Therefore, $|D|$ is base point free and hence, we have $h^1(\mathcal{O}_X(D))=0$. On the other hand, since $D.(D-C)=-3$, we have $|D-C|=\phi$ and hence, $\mathcal{O}_X(D)$ is initialized. Since $C.(C-D)=-1$, we also have $|C-D|=\phi$ and hence, we have
$$h^1(\mathcal{O}_X(C-D))=-\chi(\mathcal{O}_X(C-D))=0.$$
\noindent Moreover, since $(2C-D)^2=-2$ and $C.(2C-D)=3$, by Lemma 2.1, the movable part of $|2C-D|$ is empty. Since we have $h^0(\mathcal{O}_X(2C-D))=1$, by Riemann-Roch theorem, we have $h^1(\mathcal{O}_X(2C-D))=0$. Therefore, by Lemma 2.2, $\mathcal{O}_X(D)$ is ACM.

\smallskip

We consider the case where $D^2=4$. First of all, we show that $|D|$ is base point free. Assume that $|D|$ is not base point free, and let $\Delta$ be the base divisor of $|D|$. Let $D^{'}=D-\Delta$. Assume that ${D^{'}}^2=0$. Then, there exists an elliptic curve $F$ with $D^{'}\sim rF\;(r\geq1)$. Since $C$ is ample, we have $C.D^{'}<6$ and hence, by Lemma 2.1, we get $D^{'}\sim F$. Since $h^1(\mathcal{O}_X(D^{'}))=0$, we have 
$$\chi(\mathcal{O}_X(D^{'}))=h^0(\mathcal{O}_X(D^{'}))=h^0(\mathcal{O}_X(D))\geq\chi(\mathcal{O}_X(D)),$$
\noindent and hence, we have ${D^{'}}^2\geq D^2=4$. This is a contradiction.

Assume that ${D^{'}}^2>0$. Since $h^1(\mathcal{O}_X(D^{'}))=0$, we have ${D^{'}}^2\geq D^2=4$. Hence, by Hodge index theorem, we have $C.D^{'}\geq4$. Hence, we have $C.D^{'}=4$ or 5. If $C.D^{'}=4$, we have $(C-D^{'})^2\geq0$ and $C.(C-D^{'})=0$. Then, since $C$ is ample, we have $C\sim D^{'}$ and hence, we have $D-C\sim\Delta$. However, this contradicts to the assumption that $|D-C|=\phi$. Therefore, we have $C.D^{'}=5$. Since we have $(D^{'}-C)^2\geq-2$ and $C.(D^{'}-C)=1$, we have $|D^{'}-C|\neq\phi$. This means that $|D-C|\neq\phi$. However, by the same reason, this is a contradiction. Therefore, $|D|$ is base point free and hence, we have $h^1(\mathcal{O}_X(D))=0$. Since $(C-D)^2=-4$ and $|C-D|=|D-C|=\phi$, we have
$$h^1(\mathcal{O}_X(C-D))=-\chi(\mathcal{O}_X(C-D))=0.$$
\noindent On the other hand, $|2C-D|=\phi$. Since $C.(D-2C)=-2$, we have $|D-2C|=\phi$. Therefore, we have
$$h^1(\mathcal{O}_X(2C-D))=-\chi(\mathcal{O}_X(2C-D))=0.$$
\noindent By Lemma 2.2, $\mathcal{O}_X(D)$ is ACM.$\hfill\square$

$\;$

\noindent {\bf{Acknowledgements}}. The author would like to thank Prof. Ohashi. The author is partially supported by Grant-in-Aid for Scientific Research (25400039), Japan Society for the Promotion Science.

\end{document}